\documentclass[12pt,oneside,draft]{amsart}%
\usepackage{amsfonts, amssymb, amscd, euscript}
\usepackage[english]{babel}

\makeindex

\makeatletter \@addtoreset{equation}{section}
\@addtoreset{figure}{section} \makeatother
\newcommand{\lra}{\longrightarrow}

\newcommand{\CC}{{\mathbb{C}}}

\newcommand{\ZZ}{{\mathbb{Z}}}
\newcommand{\QQ}{{\mathbb{Q}}}

\newcommand{\NN}{{\mathbb{N}}}

\newcommand{\WWW}{{\mathcal{W}}}

\newtheorem{theorem}[subsection]{Theorem}

\newtheorem{proposition}[subsection]{Proposition}
\newtheorem{lemma}[subsection]{Lemma}
\newtheorem{corollary}[subsection]{Corollary}

\begin{document}

\title[]{Blow-ups of three-dimensional terminal singularities: $cA$ case}
\author{I. Yu. Fedorov}

\begin{abstract}
Divisors with minimal discrepancy over cA points are classified.
\end{abstract}
\maketitle  The problem of  birational classification of algebraic
varieties is highly interconnected with the problem of description
of singularities on them. One of the most important class of
three-dimensional singularities is terminal singularities, which
arise within minimal models programm. Despite the analytical
classification of the singularities
\cite{Dan82},\cite{Reid83},\cite{Ms84},\cite{Mori85}, this
description does not help one to fully understand many birationl
properties of them. In particular, the problem of description of
resolution of such singularities and the problem of classification
of morphisms of terminal varieties are still up-to-date.
Divisorial contractions to cyclic quotient singularities were
described by Y.Kawamata \cite{Kaw1}, S.Mori \cite{Mori} and S.
Cutkosky \cite{Cut} classified contractions from terminal
Gorenstein threefolds. T.Luo \cite{Luo} set out contractions when
the index is not increase. Recently M.Kawakita \cite{Kawk0}, \cite{Kawk1}, \cite{Kawk2} has gave a
description of contractions to a smooth and $cA$ points.
\par
In the paper we describe divisors with minimal discrepancy in
Mori's category over $cA$ points after M.Hayakawa \cite{Hay}, when
he did the same for the cases of non Gorenstein terminal
singularities.
\par
The author would like to thank Professor V.A. Iskovskikh and
Professor Yu.G. Prokhorov for their fruitful discussions and
encouragement. The author was partially supported by grants RFBR-99-01-01132, RFBR-96-15-96146  and
INTAS-OPEN-97-2072.

\section{\bf Preliminary results}
We will deal with varieties over $\CC$. The basic results and
notions are contained in \cite{YPG}, \cite{AVGZ}, \cite{Hay}.

\subsection{}
A singularity has $cA$ type if there is  an embedding  $j:
X\simeq\{xy+f(z,u)=0\}\hookrightarrow \CC^4_{x,y,z,u}$ (see
\cite{Mar2}). This embedding we will call  {\it standard.\/}

\subsection{}
In the paper we will study morphisms  $\pi:\bar{X}\to X$, where
$\bar{X}$ has just terminal singularities, $X$ is to be a germ of
 $cA$ point, $\pi$-exceptional set  is to be a ireducible reduced
 divisor with discrepancy 1. Such a morphisms we will call  {\it divisorial blow-ups
 with discrepancy 1.\/}

\subsection{}
Holomorphic function  $f:(\CC^n,0)\to(\CC,0)$ is to be {\it
quasi-homogeneous function\/} of degree $d$ with indices
 $\alpha_1,\dots,\alpha_n$, if for every  $\lambda>0$
we have: $f(\lambda^{\alpha_1}x_1,\dots,\lambda^{\alpha_n}x_n)=
\lambda^df(x_1,\dots,x_n)$.

The monom $\mathbf{x}^k=x_1^{k_1}\dots x_n^{k_n}$ has
{\it degree (or weight) d\/}, if $<\mathbb{
\alpha},\mathbf{k}>=\alpha_1k_1+\dots+\alpha_nk_n=d$.

Let fix the type of quasi-homogeneity $\mathbb \alpha$. Some polynomial is of the {\it order d\/}, if all of the monoms
contained in the polynomial has degree $d$ or greater.

Polynomials of order $d$ form linear space  $\mathbb{
A}_d$; $\mathbb{A}_{d'}\subset\mathbb{A}{_d}$ if $d<d'$.

Biholomorphic map $g$ is to have  order $d$, if for every
$\lambda$
\begin{equation*}
(g^*-1)\mathbb A_\lambda\subset\mathbb A_{\lambda+d}.
\end{equation*}

\subsection{}
Let $j: X\simeq\{xy+f(z,u)=0\}\hookrightarrow \CC^4_{x,y,z,u}$
be a standard embedding, the order of $f(z,u)$ is equal to $k$,
$j':X\hookrightarrow\CC^4_{x',y',z',u'}$ be some embedding $X$
such that $\chi\circ j=j'$, where $\chi$ is biholomorphic map of order  $k$ $\chi:\CC^4_{(x,y,z,u)}\to
\CC^4_{(x',y',z',u')}$. Then, by {\it pseudo weighted valuation} $\nu'$ we will mean the
pair: $\nu'=(j',\sigma')$,
where $\sigma'=(a,b,c,d)$ is weight. This pair defines weighted blow-up
$X$ which we will call $\nu'$-blow-up. If the exceptional divisor of
$\nu'$-blow-up is irreducible, then we will call $\nu'$ by {\it weighted valuation} (see \cite{Hay}).

Let's define a partial order  $\prec$ on the set of pseudo-weighted valuations.
Let $\nu'=(j',\sigma')$ and $\nu''=(j'',\sigma'')$
are two pseudo-weighted valuations, then $\nu'\succ\nu''$ if for any biholomorphic map of order $k$
$\chi:(x'',y'',z'',u'')\rightarrow (x',y',z',u')$ such that
$\chi\circ j''=j'$ and for every $f\in\CC\{x',y',z',u'\}$ the inequality $\sigma'\,wt(f)\le\sigma''\,wt(\chi^*f)$ holds.
If $\nu'\succ\nu''$ и $\nu'\prec\nu''$ then $\nu'=(j',\sigma')\sim
\nu''=(j'',\sigma'')$ is an equivalence relation.

\subsection{}
Let $\nu'=(j',\sigma')$ be a pseudo weighted valuation, where
$\sigma'=(a,b,c,d)$, $\varphi'$ is equation which defines
$j':X=\{\varphi'=0\}\hookrightarrow Y=\CC^4_{x',y',z',u'}$. Then we can construct
$\sigma'$-blow-up $\bar\pi:\bar Y\lra
Y=\CC^4_{x',y',z',u'}$ and $\nu'$-blow-up $\pi:\bar X\lra X$. Let
$\bar E$ be $\sigma'$-exceptional divisor. Then:
\begin{gather*}
K_{\bar Y}=\pi^*K_Y+(\sigma'wt(x')+\sigma'wt(y')+
\sigma'wt(z')+\sigma'wt(u')-1)\bar E \\
K_{\bar X}=\pi^*K_X+(\sigma'wt(x')+\sigma'wt(y')+
\sigma'wt(z')+\sigma'wt(u')-\\
-\sigma'wt(\varphi')-1){\bar E}_{\bar X}
\end{gather*}
Let
\begin{gather*}
d(\nu')=\sigma'wt(x')+\sigma'wt(y')+
\sigma'wt(z')+\sigma'wt(u')-\sigma'wt(\varphi')-1\\
=a+b+c+d-\sigma'wt(\varphi')-1.
\end{gather*}
We will call $d(\nu')$ by {\it virtual discrepancy
$\nu'$} (see \cite{Hay}). Note, that if ${\bar E}_{\bar X}$
is irreducible and reduced then $d(\nu')=a({\bar E}_{\bar X},X)$.

\subsection{}
For a positive rational number  $\alpha$ let's define set
\begin{equation*}
\WWW_\alpha=\{\nu'=(j',\sigma')|\  d(\nu')=\alpha\}
\end{equation*}
Let's fix an embedding  $j:X\hookrightarrow\CC^4_{x',y',z',u'}$. Then
$\WWW_\alpha(j)$ will form a subset of $\WWW_\alpha$
with embedding  $j$.

The relations  $\succ$ and $\sim$ difine partial order and equivalence relation on  $\WWW_\alpha$.

\section{\bf Main theorem}
Our main theorem is
\begin{theorem}
Let $X$ be a germ of 3-dimentional terminal $cA$ point. Then the following holds:
\begin{enumerate}
\item If $\nu\in\WWW_1$ is maximal element with respect to
$\succ$, then the $\nu$-blow-up $X$ is divisorial with discrepancy 1.
\item For every divisorial blow-up  $\pi:\bar X\to X$ with discrepancy 1,there is some
$\nu\in\WWW_1$,
such that  $\pi$ is isomorphic to the $\nu$-blow-up of $X$.
\item There is one-to one correspondence between the set of all maximal elements
$\WWW_1/\sim$ and the set of all isomorphism classes of divisorial blow-ups of $X$ with discrepancy 1.
\end{enumerate}
\end{theorem}
\begin{corollary}
Let $X$ be a germ of a 3-dimentional terminal $cA$ point, $n$ ia a number of divisors with discrepancy 1 over $X$.
Using the notation 1.1, we have that $n=deg_{min}(f)-1$, where $deg_{min}(f)$ is minimal degree among
degrees of all the monoms in
$f$.
\end{corollary}

To prove the main theorem we will
use proofs of some results from \cite{Hay} with small changes
concerning the category of singularities and equivalence relation in $W_1$.

\begin{proposition}
Let $X$ be a germ of 3-dimensional terminal Gorenstein singularity,
$\pi:\overline X\to X$ is divisorial blow-up with discrepancy 1, $E$ is $\pi$-exceptional divisor.
Let $\nu:Z\to \overline X$ be a partial resolution $\overline X$,
$\sum F_i$ is $\nu$-exceptional divisor. If
$\nu^*(E)=\nu^{-1}(E)+\sum a_iF_i$, then
\begin{equation*}
a(F_i,X)=a(F_i,\overline X)+a_i
\end{equation*}
\end{proposition}
\proof Since $K_{\bar X}=\pi^*K_X+ E$, $K_Z=\nu^*(K_{\bar
X})+\sum b_i F_i$, we have
$K_Z=\nu^*(\pi^*(K_X))+\nu^{-1}(E)+\sum (a_i+b_i)F_i$. \qed


\begin{lemma}[\cite{Hay}]
Let $\nu=(j,\sigma),\ \nu'=(j',\sigma')\in \WWW_1$ such that
$\nu\prec\nu'$. If $\nu$-blow-up $X$ is divisorial with discrepancy 1 then $\nu \sim\nu'$.
\end{lemma}
\proof is easy consequence from \cite{Hay}. \qed
\begin{lemma}
Let $\nu\sim\nu'\in \WWW_1$. Then $\nu$-blow-up of $X$ and
$\nu'$-blow-up of $X$ are isomorphic.
\end{lemma}
\proof  They are isomorphic since the blowing up ideals coincide. \qed

Let $X$ be a germ of 3-dimensional terminal $cA$ point, $j:
X\simeq\{\varphi=(xy+f(z,u))=0\}\subseteq\CC^4_{(x,y,z,u)}$ be a standard embedding,
$k$ is a minimal degree among monoms in $f$. Then the following lemma is true.
\begin{lemma}
Let $j': X\simeq\{\varphi'=0\}\subseteq\CC^4_{(x',y',z',u')}$.
If $(j',\sigma')\in \WWW_1$ then $\sigma' wt(z')=\sigma'
wt(u')=1$.
\end{lemma}
\proof Either $x'y'$, or $x'$$^2, y'$$^2\in \varphi'$ because otherwise
$j'$ will not define singularity of $cA$ type. So, $\sigma'
wt(x')+\sigma' wt(y')\geq\sigma'wt(\varphi')$. Thus we get $1=\sigma'
wt(x')+\sigma' wt(y')+\sigma' wt(z')+\sigma'
wt(u')-\sigma'wt(\varphi')-1\geq\sigma' wt(z')+\sigma' wt(u')-1$.
Hence $\sigma' wt(z')=\sigma' wt(u')=1$.\qed

Let  $\sigma_{a,b}=(a,b,1,1)$, $\nu_{a,b}=(j,\sigma_{a,b})$,
where $a,b$ are positive integers.
\begin{proposition}[\cite{Hay}]
$\WWW_1(j)=\{\nu_{a,b}|\quad a+b\leq k\}$, where $j$ is a standard embedding. In particular, $\nu_{a,b}$
is maximal in $\WWW_1(j)$ if $a+b=k$.
\end{proposition}
\proof If $(j,\sigma)\in \WWW_1(j)$, then,  $\sigma
wt(z)=\sigma wt(u)=1$ and $\sigma wt(xy)=\sigma wt(\varphi)$. Then we have
$k\geq \sigma wt(\varphi)=a+b$. \qed

\begin{theorem}
For every $\nu_{a,b}\in \WWW_1(j)$ with $a+b=k$,
the $\nu_{a,b}$-blow-up $\pi_{a,b}:\overline{X}_{a,b}\to X$ is divisorial with discrepancy 1.
These $\pi_{a,b}$ are not mutually isomorphic over $X$ and realize all possible divisors with discrepansy 1 over $X$.
\end{theorem}
\proof Lets check the terminality of $\overline{X}_{a,b}$. It is covered by four affine open charts
 $U_1,\ U_2,\ U_3,\ U_4$,
\begin{gather*}
U_1=\{\bar y+\bar x^{-k}f(\bar x\bar z, \bar x\bar u)=0\}/\ZZ_a(1,-b,-1,-1),\\
U_2=\{\bar x+\bar y^{-k}f(\bar y\bar z, \bar y\bar u)=0\}/\ZZ_b(-a,1,-1,-1),\\
U_3=\{\bar x\bar y+\bar z^{-k}f(\bar z, \bar z\bar u)=0\}\subseteq\CC^4,\\
U_4=\{\bar x\bar y+\bar u^{-k}f(\bar z\bar u, \bar
u)=0\}\subseteq\CC^4.
\end{gather*}

From this description one can see that on $\overline{X}_{a,b}$ lie just terminal singularities. In particular,
there are two cyclic factor singularities $Q_1,\ Q_2$ of
types $\frac 1 a(-1,1,1)$ and $\frac 1 b(-1,1,1)$ in charts $U_1$ and $U_2$
and at worst isolated $cDv$ points in charts $U_3$, $U_4$.

Since $E_{a,b}$ is irreducible, we have $K_{\overline
X_{a,b}}=\pi^*_{a,b}K_X+E_{a,b}$. Therefore $\pi_{a,b}$ is divisorial with discrepancy 1.

Let $D$ be the $\QQ$-Cartier divisor on $X$, defined by the equation $x=0$.
Then we have $\pi^*_{a,b}(D)=\pi^{-1}_{a,b}(D)+aE_{a,b}$.
Therefore $\pi_{a,b}$ are not mutually isomorphic over $X$.\qed

\begin{proposition}
There are just $k-1$ mutually non isomorphic divisors with discrepancy 1.
\end{proposition}
\proof Let's look on some valuation $\nu_{a,b}\in \WWW_1$,
where $a+b=k$, and blow $X$ up  $\pi_{a,b}:\overline X_{a,b}\to X$. Since  2.3,
divisors with discrepancy 1 over $X$ can not lie over points with index 1 on $\overline X_{a,b}$.
So we will study singularities with index $\geq2$. We can make so called "economic" resolution
of $Q_1$ (see \cite {YPG}), i.e. there is a projective morphism
 $\nu:Z\to \overline X_{a,b}$ such that
\begin{equation*}
K_Z=\nu^*K_{\overline K_{a,b}}+\sum_{i=1}^{a-1}\frac i a F_i,
\end{equation*}
where  $\sum F_i$  is an exceptional divisor over $Q_1$. Thus we have
\begin{equation*}
\nu^*(E_{a,b})=\nu^{-1}(E_{a,b})+\sum_{i=1}^{a-1}\frac {a-i} a
F_i.
\end{equation*}

Therefore
\begin{equation*}
a(F_i,X)=\frac i a+\frac {a-i} a=1,
\end{equation*}
for every $i=1,\dots, a-1$. On the same way, there are just $b-1$ divisors with discrepancy 1 over $Q_2$.
Thus, we have that there are only $(a-1)+(b-1)+1=k-1$ divisors over $X$ with discrepancy 1.\qed

\begin{proposition}
If $\nu'=(j',\sigma')\in \WWW_1$ is a maximal, then there is $\nu_{a,b}$, where $a+b=k$ such that $\nu'\succ\nu_{a,b}$.
\end{proposition}
\proof Let $ A=\{\nu_{a,b}\in \WWW_1|\quad
\nu_{a,b}\prec\nu'\}$. Let's check that $A$ is not empty. We have $\chi\circ j=j'$, where
 $\chi$ is some biholomorphic map with order $k$, $j$ is standard embedding. Let's look at
 weight $\sigma=(a',b',1,1)$ such that $\nu=(j,\sigma)\preceq(j',\sigma')$.
It is clear that we can choose such a weight, thus $A$ is not empty. Therefore, there is a maximal element
$\nu_{a,b}\in A$. Assuming that $a+b<k$, we will derive a contradiction.

Let $\chi$ be a quasi homogeneous map with order $k$
\begin{equation*}
\chi:\CC^4_{(x',y',z',u')}\to \CC^4_{(x,y,z,u)}
\end{equation*}
such that $\chi\circ j'=j$. Let denote
\begin{equation*}
p=\chi^*(x),\quad q=\chi^*(y),\quad r=\chi^*(z),\quad
s=\chi^*(u)\in\CC\{x',y',z',u'\}.
\end{equation*}
Then $\varphi'=pq+f(r,s)$ is the defining equation of $j'$. Since $\nu'\succ\nu_{a,b}$, we have  $\sigma'
wt(p)\geq a,\ \sigma'wt(q)\geq b$. If at least one of these inequalities are strict, then $\nu'\succ\nu_{a+1,b}$ or
$\nu'\succ\nu_{a,b+1}$ which contradict with maximality of $\nu_{a,b}$. Therefore
$\sigma' wt(p)=a,\ \sigma'wt(q)=b$. Since $a+b<k$ and  $\sigma'wt(f(s,r)\geq k$, then
$\sigma'wt(\varphi')=a+b$. It follows from 2.6 that up to permutation of coordinates
\begin{equation*}
\sigma'wt(x',y',z',u')=(a',b',1,1)
\end{equation*}
for some $a',b'\in\NN$. Since $d(\nu')=1$, we have
$a'+b'=a+b$. It follows from 2.4, that after a permutation $x'$ and
$y'$, we get $\nu'\sim\nu_{a,b}$ which also contradicts the maximality of
$\nu'$. Hence $a+b=k$. \qed

 The main theorem is proved since, from 2.4 and 2.9 we see that  $\nu_{a,b}$, where $a+b=k$ are all $k-1$
 maximal elements in $\WWW_1$. On the other hand, from 2.7 and 2.8 follows that over
 $X$ there are only $k-1$ divisor with discrepancy 1.

\end{document}